\documentclass{amsart}
  \usepackage{amsmath}
\usepackage{amssymb}
\usepackage{amsfonts}
 \usepackage{eucal}
    \usepackage[monochrome]{xcolor}
 \makeatletter
     \def\section{\@startsection{section}{1}%
     \z@{.7\linespacing\@plus\linespacing}{.5\linespacing}%
     {\bfseries
     \centering
     }}
     \def\@secnumfont{\bfseries}
     \makeatother
   \newtheorem{theorem}{Theorem}[section]

\newtheorem{corollary}[theorem]{Corollary}

\theoremstyle{definition}

\newtheorem{remark}[theorem]{Remark}

\numberwithin{equation}{section}

\def \D{{\Delta}}

\def \e{{\varepsilon}}

\def \l{{\lambda}}

\def \s{{\sigma}}

\def \qq{{\qquad}}

\def\beq{\begin{equation}}
\def\eeq{\end{equation}}
  \def\ben{\begin{eqnarray}}
\def\een{\end{eqnarray}}

\def\P{{\mathbb P}}

\def\R{{\mathbb R}}

\def\Z{{\mathbb Z}}

 \scrollmode
  \font\sevenrm= cmr10 at 7 pt

\def\ddate {\sevenrm \ifcase\month\or January\or
February\or March\or April\or May\or June\or July\or
August\or September\or October\or November\or December\fi\! {\the\day}, \!{\sevenrm\the\year}}

\date{\today}

 \title[\rm 
  On Mukhin's NSC for the local limit theorem]
  {On Mukhin's necessary and sufficient condition for the validity of the local limit theorem}

 \author{Michel  J.\,G. WEBER}
\address{IRMA, UMR 7501, Universit\'e
Louis-Pasteur et C.N.R.S.,   7  rue Ren\'e Descartes, 67084
Strasbourg Cedex, France.
   E-mail:    {\tt  michel.weber@math.unistra.fr}}

\begin{document}
 \maketitle


\renewcommand{\thefootnote}{} {{
 \renewcommand{\thefootnote}{\arabic{footnote}}
\setcounter{footnote}{0}
  \begin{abstract}   Mukhin found in 1984 an important necessary and sufficient condition for the validity of the local limit theorem.  Revisiting the succint proof given in   \cite{Mu2}, we could only prove rigorously a weaker necessary and sufficient condition,   with a significantly different formulation. This is  the object of this short Note.
   \end{abstract}

\section{Introduction-Result.}\label{s1}

  \maketitle 
     Let $\{S_n,n\ge 1\}$ be a sequence of  integer-valued random variables  such  that an integral limit theorem  holds:  there exist $a_n\in \R$ and real $b_n\to \infty$ such that the sequence of distributions of $(S_n-a_n)/b_n$ converges weakly to an absolutely continuous distribution $G$ with density $g(x)$, which is uniformly continuous in $\R$\,\footnote{It was assumed in \cite{Mu} but not in \cite{Mu2} that $g$ is uniformly continuous in $\R$.}.  
The  local limit theorem is   valid if
\beq\label{ilt.llt}
\P\{S_n=m\}=B_n^{-1} g\Big(\frac{m-A_n}{B_n}\Big) + o(B_n^{-1}),
\eeq
uniformly in $m\in \Z$.

It is known that the local limit theorem always implies the integral limit theorem. The converse is false as was proved by Gamkrelidze \cite{Gam2a}, who showed  that there exists a sequence of independent centered, square integrable, integer-valued random variables, with partial sums $S_n$, $\s_n^2={\rm Var}(S_n)$, that satisfy an integral limit theorem,  are asymptotically uniformly distributed: $\P\{S_n \equiv m \,{\rm (mod)}\,d\}\to 1/d$, as $n\to \infty$,
for any  $m = 0,1,\ldots,d-1$, $d\ge 2$,   and further  $S_n/\s_n$ asymptotically is normal, but   the local limit theorem does not hold. By Rozanov's theorem, the asymptotic uniform distribution property is a necessary condition for the validity of the local limit theorem, see \cite{SW}. This is recently  improved in  \cite{W}.

\vskip 5 pt  Mukhin \cite[Th.\,1]{Mu2} has shown  the following important result   relating the integral limit theorem to  the  local limit theorem.

\begin{theorem}\label{Mukhin[NSC]}   The following assertions are equivalent.
\vskip 3 pt {\rm (A)} \quad There exists a sequence of integers $v_n=o(b_n)$ such that
\beq \sup_{m}\Big|\P\{S_n=m+v_n\big\}-\P\{S_n=m \big\}\Big|\,=\,\,o\Big(\frac{1}{b_n}\Big)
\eeq
\vskip 3 pt
{\rm (B)} \quad 
\beq\P\{S_n=m\}\,=\, \frac{1}{b_n} \,g\Big(\frac{m-a_n}{b_n}\Big) +  \,o\Big(\frac{1}{b_n}\Big),
\eeq
 \end{theorem}
 
 \vskip 10 pt  However from the proof available   in  \cite{Mu2},  we could only  obtain the following significantly different and weaker result.  By the integral limit theorem we have 
\beq \label{en}\e_n:=\sup_{x\in \R}\Big|\P\Big\{\frac{S_n-a_n}{b_n}<x\Big\}-G(x)\Big|\ \to \ 0.
\eeq

\begin{theorem}\label{Mukhin[NSC]weaker}
Let $v_n$ be a  sequence of positive  integers such that $v_n=o(b_n)$.  The following assertions are equivalent.
\vskip 3 pt
 {\rm (A')} \quad  \beq \sup_{m, k\in\Z\atop |m-k|\le v_n}\Big|\P\{S_n=m\big\}-\P\{S_n=k \big\}\Big|\,=\,\,o\Big(\frac{1}{b_n}\Big)+ {\color{blue} \frac{1}{v_n}\,\mathcal O(\e_n)},\qq \quad n\to \infty.
\eeq

\vskip 3 pt
{\rm (B')} \quad 
\beq\sup_{m}\Big|\P\{S_n=m\} -\frac{1}{b_n} \,g\Big(\frac{m-a_n}{b_n}\Big)\Big| \,=\,   \,o\Big(\frac{1}{b_n}\Big)+ {\color{blue} \frac{1}{v_n}\,\mathcal O(\e_n)},\qq \quad n\to \infty.\eeq
 \end{theorem}
\vskip 3 pt
Choosing 
$v_n= \max\{1, [\sqrt \e_n b_n]\}$ 
 we get the following

\begin{corollary}\label{Mukhin[NSC]w}  A necessary and sufficient condition for the local limit theorem in the usual form to hold is 
\beq \sup_{m, k\in\Z\atop |m-k|\le \max\{1, [\sqrt \e_n b_n]\}}\Big|\P\{S_n=m\big\}-\P\{S_n=k \big\}\Big|\,=\,\,o\Big(\frac{1}{b_n}\Big).
\eeq
 \end{corollary}
\vskip 2 pt 
 \begin{remark} 
By Theorem 2 in \cite{Mu2}, for arbitrary $\D>0$, the relation 
\beq \label{(3)}\P\{S_n\in[x,x+\D]\}= \frac{\D}{b_n}\,g\Big(\frac{m-a_n}{b_n}\Big) +  \,o\Big(\frac{1}{b_n}\Big),\qq \quad n\to \infty,
\eeq
holds uniformly in $x\in\R$ if and only if for any $\l>0$, $v_n>0$, $v_n=o(b_n)$ we have
\beq  \label{(4)}\sup_{x } \Big|\P\big\{S_n\in[x+v_n,x+v_n+\l]\big\}-\P\big\{S_n\in[x+v_n,x+v_n]\big\}\Big|=  \,o\Big(\frac{1}{b_n}\Big),\quad\,    n\to \infty.
\eeq
It is not clear how these assertions should imply our condition (A'), and conversely.
\end{remark}
\vskip 3 pt 

\begin{remark}  
   The key quantity 
\beq \sup_{m \in\Z }\big|\P\{S_n=m+k\big\}-\P\{S_n=k \big\}\big| ,
\eeq for partial sums of integer valued random variables was thoroughly investigated by Mukhin in several works,   in \cite{Mu3}, \cite{Mu1}, \cite{Mu} and \cite{Mu2a} notably, using structural characteristics of the summands. The case  of analogous differences for densities and arbitrary distributions was also considered. However Mukhin wrote in \cite{Mu2}: \lq\lq ... getting from here more general sufficient conditions turns out to be difficult in view of the lack of good criteria for relations \eqref{(3)} and \eqref{(4)}.  Working with asymptotic equidistribution properties are more convenient in this respect\,\rq\rq.  Sufficient conditions of different type for the validity of the  local limit theorem are given in \cite{Mu2b}. In To An'\,Zung \cite{VMT1},    some upper bounds for $\sup_{j\in \Z^s} \big|\P\{S_n= j+r\}-\P\{S_n= j\}\big| $ are  proved,  where $S_n= X_1+\ldots +X_n$, $X_n$ being independent $\Z^s$-valued  random vectors, and some local limit theorems are also presented.
\end{remark}
\section{Proof.} (A')$\Rightarrow$(B'). {\color{blue}Let $m$ be arbitrary.} We write
 \ben v_n\P\{S_n=m\}&=&\sum_{k=m}^{m+v_n-1} \P\{S_n=k\}+ \sum_{k=m}^{m+v_n-1}\big(\P\{S_n=m\}-\P\{S_n=k\}\big) 
 \cr &=& {\rm (I)} +{\rm (II)}.
 \een
We have
\ben  {\rm (I)} &=& \P\{ m\le S_n\le m+v_n-1\}
\cr &=& G\Big(\frac{m+v_n-1-a_n}{b_n}\Big)-
G\Big(\frac{m-a_n}{b_n}\Big)+ \mathcal O(\e_n)
\cr &=& \frac{v_n-1}{b_n} \,g\Big(\frac{x_{n,m}-a_n}{b_n}\Big)+  \mathcal O(\e_n),
 \een
for some $x_{n,m}\in ]m ,m+v_n-1[$.  As $g$   is uniformly continuous in $\R$, 
{\color{blue}$$ \forall \e>0, \exists \eta>0\,:\,
\ \sup_{x,y\in\R\atop |x-y|\le \eta}|g(x)-g(y)|\le \e,$$   
 and $v_n=o(b_n)$,  it follows that given any $\e>0$, if $n$ is large enough so that 
$  \frac{v_n}{b_n}\le \eta$, then 
$$\sup_{u,v\in\R\atop |u-v|\le {v_n} }\Big|g\big(\frac{u}{b_n}\big)-g\big(\frac{v}{b_n}\big)\Big|\le \e.$$
Therefore  
\beq \label{unif.g}\sup_{u,v\in\R\atop |u-v|\le {v_n} }\Big|g\big(\frac{u}{b_n}\big)-g\big(\frac{v}{b_n}\big)\Big|=o_g(1), \qq \quad n\to \infty.
\eeq
Whence also 
$$\sup_{m, n\in\Z\atop |m-n|\le v_n}\Big|g\big(\frac{x_{n,m}-a_n}{b_n}\big)-g\big(\frac{m-a_n}{b_n}\big)\Big| =   o_g(1), \qq \quad n\to \infty.$$
}

 Thus
\ben  {\rm (I)} &=&   \frac{v_n-1}{b_n} \, g\Big(\frac{m-a_n}{b_n}\Big)\, +\, \frac{v_n-1}{b_n} \,o_g(1)+  \mathcal O(\e_n).
 \cr  |{\rm (II)}| &=& v_n\, o\Big(\frac{1}{b_n}\Big)+ \,\mathcal O(\e_n)  \qq\quad \hbox{by (A')}.
 \een By combining both estimates, next dividing by $v_n$, we get

\beq
\P\{S_n=m\}\,=\, \frac{1}{b_n} \,g\Big(\frac{m-a_n}{b_n}\Big) +  \,o\Big(\frac{1}{b_n}\Big)+ {\color{blue} \frac{1}{v_n}\,\mathcal O(\e_n)}.
\eeq

We rewrite this as follows
\beq
\Big|\P\{S_n=m\}- \frac{1}{b_n} \,g\Big(\frac{m-a_n}{b_n}\Big)\Big|\,=\,    \,o\Big(\frac{1}{b_n}\Big)+ {\color{blue} \frac{1}{v_n}\,\mathcal O(\e_n)},
\eeq
and note that this is true for all $m$, as we started the proof with $m$ arbitrary, and the analysis   made only involved   the interval $[m,m+v_n]$.



\vskip 10 pt 
(B')$\Rightarrow$(A').
Let $v_n=o(b_n)$. We have for $k,m\in\Z$ such that $|k-m|\le v_n$, using \eqref{unif.g},
\ben \P\{S_n=m \big\}-\P\{S_n=k \big\}&=&\frac{1}{b_n} \,\Big\{g\Big(\frac{m -a_n}{b_n}\Big)- g\Big(\frac{k-a_n}{b_n}\Big)\Big\}
\cr & & +       \,o\Big(\frac{1}{b_n}\Big)+\frac{1}{v_n}\,\mathcal O(\e_n) 
\cr
&=&\frac{1}{b_n} \, o_g(1)+   \,o\Big(\frac{1}{b_n}\Big)  +  \frac{1}{v_n}\,\mathcal O(\e_n)
\cr &=&    \,o\Big(\frac{1}{b_n}\Big)+  \frac{1}{v_n}\,\mathcal O(\e_n).
\een

%
 \end{document}